\documentclass[a4paper,12pt,twoside,leqno]{article} 
\usepackage[english]{babel} 
\usepackage{amsmath,amssymb,amsthm,amsfonts} 
\usepackage[dvips]{graphicx}
\usepackage[all]{xy}

\usepackage{hyperref}
\hypersetup{
  pdftitle={Steinitz classes of tamely ramified nonabelian extensions of odd prime power order},
  pdfauthor={Alessandro Cobbe},
  pdfsubject={Steinitz classes},
  pdfkeywords={Steinitz classes}}

\usepackage{fancyhdr}
\newcommand{\Q}{\mathbb Q}
\newcommand{\N}{\mathbb N}
\newcommand{\Z}{\mathbb Z}

\newcommand{\p}{\mathfrak p}

\newcommand{\gal}{\mathrm{Gal}}
\renewcommand{\epsilon}{\varepsilon}
\newcommand{\disc}{\mathrm{d}}
\newcommand{\norm}{\mathrm{N}}

\newcommand{\st}{\mathrm{st}}

\newcommand{\rt}{\mathrm{R}_t}
\newcommand{\cl}{\mathrm{Cl}}

\newcommand{\G}{\mathcal G}
\newcommand{\oo}{\mathcal O}

\newcommand{\comment}[1]{}

\newtheorem{teo}{Theorem}[section]
\newtheorem{lemma}[teo]{Lemma}

\newtheorem{prop}[teo]{Proposition}
\newtheorem{defn}[teo]{Definition}

\title{Steinitz classes of tamely ramified nonabelian extensions of odd prime power order}
\author{Alessandro Cobbe}

\begin{document}
\maketitle

\markboth{Abstract}{Abstract}
\begin{section}*{Abstract}
\addcontentsline{toc}{section}{Abstract}
The Steinitz class of a number field extension $K/k$ is an ideal class in the ring of integers $\oo_k$ of $k$, which, together with the degree $[K:k]$ of the extension determines the $\oo_k$-module structure of $\oo_K$. We call $\rt(k,G)$ the classes which are Steinitz classes of a tamely ramified $G$-extension of $k$. We will say that those classes are realizable for the group $G$; it is conjectured that the set of realizable classes is always a group.

In this paper we will develop some of the ideas contained in \cite{Cobbe} to study some $l$-groups, where $l$ is an odd prime number. In particular, together with \cite{Bruche} we will complete the study of realizable Steinitz classes for groups of order $l^3$. We will also give an alternative proof of the results of \cite{Bruche}, based on class field theory.
\end{section}

\markboth{Introduction}{Introduction}
\begin{section}*{Introduction}
\addcontentsline{toc}{section}{Introduction}

Let $K/k$ be an extension of number fields and let $\oo_K$ and $\oo_k$ be their rings of integers. By Theorem 1.13 in \cite{Narkiewicz} we know that
\[\oo_K\cong \oo_k^{[K:k]-1}\oplus I\]
where $I$ is an ideal of $\oo_k$. By Theorem 1.14 in \cite{Narkiewicz} the $\oo_k$-module structure of $\oo_K$ is determined by $[K:k]$ and the ideal class of $I$. This class is called the \emph{Steinitz class} of $K/k$ and we will indicate it by $\st(K/k)$. Let $k$ be a number field and $G$ a finite group, then we define:
\[\rt(k,G)=\{x\in\cl(k):\ \exists K/k\text{ tame, }\gal(K/k)\cong G, \st(K/k)=x\}.\]

In this paper we will use the notations and some techniques from \cite{Cobbe} to study the realizable classes for some $l$-groups, where $l$ is an odd prime number.

Some of the results in this paper are parts of the author's PhD thesis \cite{tesi}. 
For earlier results see \cite{Bruche}, \cite{BrucheSodaigui}, \cite{ByottGreitherSodaigui}, \cite{Carter}, \cite{Carter2}, \cite{CarterSodaigui_quaternionigeneralizzati}, \cite{Endo}, \cite{GodinSodaigui_A4}, \cite{GodinSodaigui_ottaedri}, \cite{Long2}, \cite{Long0},  \cite{Massy}, \cite{McCulloh}, \cite{McCulloh_Crelle}  \cite{Sodaigui1}, \cite{Sodaigui2} and \cite{Soverchia}.

\end{section}

\markboth{Acknowledgements}{Acknowledgements}
\begin{section}*{Acknowledgements}
\addcontentsline{toc}{section}{Acknowledgements}
I am very grateful to Professor Cornelius Greither and to Professor Roberto Dvornicich for their advice and for the patience they showed, assisting me in the writing of my PhD thesis with a lot of suggestions. I also wish to thank the Scuola Normale Superiore of Pisa, for its role in my mathematical education and for its support during the time I was working on my PhD thesis.
\end{section}

\begin{section}{Preliminary results}
We start recalling the following two fundamental results.
\begin{teo}\label{discriminante}
If $K/k$ is a finite tame Galois extension then 
\[\disc(K/k)=\prod_{\p}\p^{(e_\p-1)\frac{[K:k]}{e_\p}},\]
where $e_\p$ is the ramification index of $\p$.
\end{teo}

\begin{proof}
This follows by Propositions 8 and 14 of chapter III of \cite{Lang}.
\end{proof}

\begin{teo}\label{stdisc}
Assume $K$ is a finite Galois extension of a number field $k$.
\begin{enumerate}
\item[(a)] If its Galois group either has odd order or has a noncyclic $2$-Sylow subgroup then $\disc(K/k)$ is the square of an ideal and this ideal represents the Steinitz class of the extension.
\item[(b)] If its Galois group is of even order with a cyclic $2$-Sylow subgroup and $\alpha$ is any element of $k$ whose square root generates the quadratic subextension of $K/k$ then $\disc(K/k)/\alpha$ is the square of a fractional ideal and this ideal represents the Steinitz class of the extension.
\end{enumerate}
\end{teo}

\begin{proof}
This is a corollary of Theorem I.1.1 in \cite{Endo}. In particular it is shown in \cite{Endo} that in case (b) $K/k$ does have exactly one quadratic subextension.
\end{proof}

Further, considering Steinitz classes in towers of extensions, we will need the following proposition.

\begin{prop}\label{stintermediateextension}
Suppose $K/k_1$ and $k_1/k$ are number fields extensions. Then
\[\st(K/k)=\st(k_1/k)^{[k_1:E]}\norm_{k_1/k}(\st(K/k_1)).\]
\end{prop}

\begin{proof}
This is Proposition I.1.2 in \cite{Endo}.
\end{proof}

We will also use some other preliminary results.

\begin{lemma}\label{congruenza}
Let $m,n,x,y$ be integers. If $x\equiv y\pmod{m}$ and any prime $q$ dividing $n$ divides also $m$ then
\[x^n\equiv y^n\pmod{mn}.\]
\end{lemma}

\begin{proof}
Let $n=q_1\dots q_r$ be the prime decomposition of $n$ ($q_i$ and $q_j$ with $i\neq j$ are allowed to be equal). We prove by induction on $r$ that $x^n\equiv y^n\pmod{mn}$. If $r=1$, then $mn=mq_1$ must divide $m^{q_1}$ and there exists $b\in\N$ such that
\[x^n=(y+bm)^{q_1}=y^{q_1}+\sum_{i=1}^{q_1-1}\binom{q_1}{i}(bm)^i y^{q_1-i}+(bm)^{q_1}\equiv y^n\pmod{mn}.\]
Let us assume that the lemma is true for $r-1$ and prove it for $r$. Since $q_r|m$, as above, for some $c\in \N$ we have
\[\begin{split}x^n&=(y^{q_1\dots q_{r-1}}+cmq_1\dots q_{r-1})^{q_r}\\&=y^n+\sum_{i=1}^{q_r}\binom{q_r}{i}(cmq_1\dots q_{r-1})^i y^{q_1\dots q_{r-1}(q_r-i)}\equiv y^n\pmod{mn}.\end{split}\]
\end{proof}

\begin{defn}
Let $K/k$ be a finite abelian extension of number fields. Then we define the subgroup $W(k,K)$ of the ideal class group of $k$ in the following equivalent ways (the equivalence is shown in \cite{Cobbe}, Proposition 1.10):
\[\begin{split}
&W(k,K)=\{x\in J_k/P_k:\text{$x$ contains infinitely many primes of absolute}\\&\qquad\qquad \qquad \text{degree $1$ splitting completely in $K$}\}\\
&W(k,K)=\{x\in J_k/P_k:\text{$x$ contains a prime splitting completely in $K$}\}\\
&W(k,K)=\norm_{K/k}(J_K)\cdot P_k/P_k.\end{split}\]
\end{defn}
In the case of cyclotomic extensions we will also use the shorter notation $W(k,m)=W(k,k(\zeta_m))$.

\begin{lemma}\label{Wexp}
If $q|n\Rightarrow q|m$ then
$W(k,m)^n\subseteq W(k,mn)$.
\end{lemma}

\begin{proof}
Let $x\in W(k,m)$. According to Proposition 1.10 and Lemma 1.11, both from \cite{Cobbe}, $x$ contains a prime ideal $\p$, prime to $mn$ and such that $\norm_{k/\Q}(\p)\in P_\Q^\mathfrak m$, where $\mathfrak m=m\cdot p_\infty$. Then by Lemma \ref{congruenza}, $\norm_{k/\Q}(\p^n)\in P_\Q^\mathfrak n$, with $\mathfrak n=mn\cdot p_\infty$, and it follows from Lemma 1.12 of \cite{Cobbe} that $x^n\in W(k,mn)$.
\end{proof}

\begin{defn}
We will call a finite group $G$ of order $m$ \emph{good} if the following properties are verified:
\begin{enumerate}
\item For any number field $k$, $\rt(k,G)$ is a group.
\item For any tame $G$-extension $K/k$ of number fields there exists an element $\alpha_{K/k}\in k$ such that:
\begin{enumerate}
\item[(a)] If $G$ is of even order with a cyclic $2$-Sylow subgroup, then a square root of $\alpha_{K/k}$ generates the quadratic subextension of $K/k$; if $G$ either has odd order or has a noncyclic $2$-Sylow subgroup, then $\alpha_{K/k}=1$.
\item[(b)] For any prime $\p$, with ramification index $e_\p$ in $K/k$, the ideal class of
\[\left(\p^{(e_\p-1)\frac{m}{e_\p}-v_\p(\alpha_{K/k})}\right)^\frac{1}{2}\]
is in $\rt(k,G)$. 
\end{enumerate}
\item For any tame $G$-extension $K/k$ of number fields, for any prime ideal $\p$ of $k$ and any rational prime $l$ dividing its ramification index $e_\p$, the class of the ideal
\[\p^{(l-1)\frac{m}{e_\p(l)}}\]
is in $\rt(k,G)$ and, if $2$ divides $(l-1)\frac{m}{e_\p(l)}$, the class of
\[\p^{\frac{l-1}{2}\frac{m}{e_\p(l)}}\]
is in $\rt(k,G)$.
\item $G$ is such that for any number field $k$, for any class $x\in\rt(k,G)$ and any integer $a$, there exists a tame $G$-extension $K$ with Steinitz class $x$ and such that every non trivial subextension of $K/k$ is ramified at some primes which are unramified in $k(\zeta_{a})/k$.
\end{enumerate}
\end{defn}
\end{section} 

\begin{section}{Some $l$-groups}
In \cite{Bruche}, Cl\'ement Bruche proved that if $G$ is a nonabelian group of order $l^3=uv$ and exponent $v$, where $l$ is an odd prime, then $\rt(k,G)=W(k,l)^{u(l-1)/2}$ under the hypothesis that the extension $k(\zeta_v)/k(\zeta_l)$ is unramified, thereby giving an unconditional result when $G$ has exponent $l$.

In this section we prove that $\rt(k,C(l^2)\rtimes_\mu C(l))=W(k,l)^{l(l-1)/2}$, without any additional hypothesis on the number field $k$. Indeed we will consider a more general situation, studying groups of the form $G=C(l^n)\rtimes_\mu C(l)$, with $n\geq 2$, where $\mu$ sends a generator of $C(l)$ to the elevation to the $l^{n-1}+1$-th power. Together with Bruche's result this will conclude the study of realizable Steinitz classes for tame Galois extensions of degree $l^3$.

\begin{lemma}\label{identificopgruppo}
Let $l$ be an odd prime. The group $G=C(l^n)\rtimes_\mu C(l)$, with $n\geq 2$ is identified by the exact sequence
\[1\to\ C(l^n)\to G\to C(l)\to 1\]
if the action of $C(l)$ on $C(l^n)$ is given by $\mu$.
\end{lemma}

\begin{proof}
Let $G$ be the group written in the above exact sequence, let $H$ be a subgroup of $G$ isomorphic to $C(l^n)$ and generated by $\tau$; let $x\in G$ be such that its class modulo $H$ generates $G/H$, which is cyclic of order $l$, and such that $x\tau x^{-1}=\tau^{l^{n-1}+1}$, i.e. $x\tau=\tau^{l^{n-1}+1}x$. Then $x^l=\tau^a$ for some $a\in \N$. 
Since $G$ is of order $l^{n+1}$ and it is not cyclic, the order of $x$ must divide $l^n$ and so
\[\tau^{al^{n-1}}=x^{l^n}=1,\]
i.e. $l$ divides $a$ and there exists $b\in \N$ such that $a=b l$. By induction we prove that, for $m\geq 1$,
\[(\tau^{-b}x)^m=\tau^{-bm-bl^{n-1}(m-1)m/2}x^m.\]
This is obvious for $m=1$; we have to prove the inductive step:
\[\begin{split}(\tau^{-b}x)^m&=\tau^{-b(m-1)-bl^{n-1}(m-2)(m-1)/2}x^{m-1}\tau^{-b}x\\
&=\tau^{-b(m-1)-bl^{n-1}(m-2)(m-1)/2}x^{m-1}\tau^{-b}x^{-(m-1)}x^m\\
&=\tau^{{-b}(m-1)-bl^{n-1}(m-2)(m-1)/2}\tau^{-b(1+l^{n-1})^{m-1}}x^m\\
&=\tau^{-b(m-1)-bl^{n-1}(m-2)(m-1)/2-b-b(m-1)l^{n-1}}x^m\\
&=\tau^{-bm-bl^{n-1}(m-1)m/2}x^m.\end{split}\]
Then calling $\sigma=\tau^{-b}x$, we obtain that
\[\sigma^l=(\tau^{-b}x)^l=\tau^{-b l}x^l=\tau^{-a+a}=1.\]
Further
\[\sigma\tau\sigma^{-1}=\tau^{-b}x\tau x^{-1}\tau^b=\tau^{-b}\tau^{l^{n-1}+1}\tau^b=\tau^{l^{n-1}+1}\]
and $\sigma,\tau$ are generators of $G$. Thus $G$ must be a quotient of the group
\[\langle\sigma,\tau:\ \sigma^l=\tau^{l^n}=1,\sigma\tau\sigma^{-1}=\tau^{l^{n-1}+1}\rangle.\]
But this group has the same order of $G$ and thus they must be equal.
\end{proof}

It follows that we can use Proposition 2.13 of \cite{Cobbe} to study $\rt(k,C(l^n)\rtimes_{\mu} C(l))$, for any number field $k$.

For any $\tau\in H$ we define $E_{k,\mu,\tau}$ as the fixed field in $k(\zeta_{o(\tau)})$ of
\[G_{k,\mu,\tau}=\left\{g\in\gal(k(\zeta_{o(\tau)})/k):\ \exists g_1\in\G,\ \mu(g_1)(\tau)=\tau^{\nu_{k,\tau}(g_2)}\right\}\]
where $g_2(\zeta_{o(\tau)})=\zeta_{o(\tau)}^{\nu_{k,\tau}(g_2)}$ for any $g_2\in \gal(k(\zeta_{o(\tau)})/k)$.

\begin{lemma}\label{Eciclo}
Let $\tau$ be a generator of $C(l^n)$ in $C(l^n)\rtimes_{\mu} C(l)$. Then $E_{k,\mu,\tau}=k(\zeta_{l^{n-1}})$.
\end{lemma}

\begin{proof}
By definition $E_{k,\mu,\tau}$ is the fixed field in $k(\zeta_{l^n})$ of
\[\begin{split}G_{k,\mu,\tau}&=\left\{g\in\gal(k(\zeta_{l^n})/k):\ \exists g_1\in C(l)\ \mu(g_1)(\tau)=\tau^{\nu_{k,\tau}(g)}\right\}\\
&=\left\{g\in\gal(k(\zeta_{l^n})/k):\ \exists a\in \N\ \tau^{al^{n-1}+1}=\tau^{\nu_{k,\tau}(g)}\right\}\\
&=\left\{g\in\gal(k(\zeta_{l^n})/k): \nu_{k,\tau}(g)\equiv1\pmod{l^{n-1}}\right\}
\\&=\left\{g\in\gal(k(\zeta_{l^n})/k): g(\zeta_{l^{n-1}})=\zeta_{l^{n-1}}\right\}
=\gal(k(\zeta_{l^n})/k(\zeta_{l^{n-1}})).
\end{split}\]
Hence $E_{k,\mu,\tau}=k(\zeta_{l^{n-1}})$.
\end{proof}

\begin{lemma}\label{lncostruttivo}
We have
\[\rt(k,C(l^n)\rtimes_{\mu} C(l))\supseteq W(k,l^{n-1})^{\frac{l-1}{2}l}.\]
Further, for any $x\in W(k,l^{n-1})$ and any positive integer $a$, there exists a tame $G$-extension of $k$ with Steinitz class $x^{\frac{l-1}{2}l}$ and such that any nontrivial subextension of $K/k$ is ramified at some primes which are unramified in $k(\zeta_a)/k$.
\end{lemma}
\begin{proof}
By Proposition 2.13 of \cite{Cobbe} and Lemma \ref{identificopgruppo},
\[\rt(k,C(l^n)\rtimes_\mu C(l))\supseteq\rt(k,C(l))^{l^{n}} \cdot W\left(k,E_{k,\mu,\tau}\right)^{\frac{l-1}{2}l},\]
where $\tau$ is a generator of $C(l^n)$.
We easily conclude since $1\in \rt(k,C(l))$ and, by Lemma \ref{Eciclo}, $E_{k,\mu,\tau}=k(\zeta_{l^{n-1}})$, i.e. 
\[W(k,E_{k,\mu,\tau})=W(k,l^{n-1}).\]
The second part of the lemma follows again by Proposition 2.13 of \cite{Cobbe}.
\end{proof}
To prove the opposite inclusion we need some lemmas.

\begin{lemma}\label{potenzagalois}
Let $\tau$ be a generator of $C(l^n)$ in $C(l^n)\rtimes_{\mu} C(l)$ and $0<c<n$ be an integer, then
\[G_{k,\mu,\tau^{l^c}}^{l^c}\subseteq G_{k,\mu,\tau}.\]
\end{lemma}

\begin{proof}
For any positive integer $a$ we define
\[\hat\mu_{\tau^{a}}:\G\to(\Z/o(\tau^a)\Z)^*\]
by $\tau^{a\hat\mu_\tau(g_1)}=\mu(g_1)\left(\tau^{a}\right)$ for all $g_1\in\G$. By definition, if $g\in G_{k,\mu,\tau^c}$, then there exists $g_1\in\G$ such that \[\tau^{l^c\nu_{k,\tau^{l^c}}(g)}=\mu(g_1)(\tau^{l^c})=\tau^{l^c\hat\mu_{\tau^{l^c}}(g_1)}.\]
We also observe that
\[\zeta_{l^{n-c}}^{\nu_{k,\tau}(g)}=\zeta_{l^n}^{l^c\nu_{k,\tau}(g)}=g\left(\zeta_{l^n}\right)^{l^c}=g\left(\zeta_{l^{n-c}}\right)=\zeta_{l^{n-c}}^{\nu_{k,\tau^{l^c}}(g)}\]
and that
\[\tau^{l^c\hat\mu_{\tau^{l^c}}(g_1)}=\mu(g_1)(\tau^{l^c})=\mu(g_1)(\tau)^{l^c}=\tau^{l^c\hat\mu_\tau(g_1)}.\]
From the above equalities we deduce that
\[\nu_{k,\tau}(g)\equiv\nu_{k,\tau^{l^c}}(g)\equiv \hat\mu_{\tau^{l^c}}(g_1)\equiv\hat\mu_{\tau}(g_1)\pmod{l^{n-c}}\]
and therefore by Lemma \ref{congruenza} we obtain that
\[\nu_{k,\tau}(g^{l^c})\equiv\hat\mu_{\tau}(g_1^{l^c})\pmod{l^n}.\]
We conclude that
\[\tau^{\nu_{k,\tau}(g^{l^c})}=\tau^{\hat\mu_{\tau}(g_1^{l^c})}=\mu(g_1^{l^c})(\tau)\]
and hence that $g^{l^c}\in G_{k,\mu,\tau}$.
\end{proof}

\begin{lemma}\label{troppo}
Let $\tau$ be a generator of $C(l^n)$ in $C(l^n)\rtimes_{\mu} C(l)$ and $0<c<n$ be an integer, then
\[W(k,E_{k,\mu,\tau^{l^c}})^{l^c}\subseteq W(k,l^{n-1}).\]
\end{lemma}

\begin{proof}
Let $x$ be a class in $W(k,E_{k,\mu,\tau^{l^c}})$. By Proposition 1.10 in \cite{Cobbe} there exists a prime $\p$ in the class of $x$ splitting completely in $E_{k,\mu,\tau^{l^c}}/k$. By Theorem IV.8.4 in \cite{Neukirch}, $\p\in H_{E_{k,\mu,\tau^{l^c}}/k}^\mathfrak m$, where $\mathfrak m$ is a cycle of declaration of $E_{k,\mu,\tau^{l^c}}/k$.
Then, by Proposition II.3.3 in \cite{Neukirch},
\[\left.\left(\frac{k(\zeta_{l^n})/k}{\p}\right)\right|_{E_{k,\mu,\tau^{l^c}}}=\left(\frac{E_{k,\mu,\tau^{l^c}}/k}{\p}\right)=1.\]
Thus
\[\left(\frac{k(\zeta_{l^n})/k}{\p}\right)\in\gal(k(\zeta_{l^n})/E_{k,\mu,\tau^{l^c}})=G_{k,\mu,\tau^{l^c}}\]
and it follows by Lemma \ref{potenzagalois} that
\[\left(\frac{k(\zeta_{l^n})/k}{\p^{l^c}}\right)=\left(\frac{k(\zeta_{l^n})/k}{\p}\right)^{l^c}\in G_{k,\mu,\tau^{l^c}}^{l^c}\subseteq G_{k,\mu,\tau}=\gal(k(\zeta_{l^n})/E_{k,\mu,\tau}).\]
Then
\[\left(\frac{E_{k,\mu,\tau}/k}{\p^{l^c}}\right)=\left.\left(\frac{k(\zeta_{l^n})/k}{\p^{l^c}}\right)\right|_{E_{k,\mu,\tau}}=1\]
and so the class $x^{l^c}$ of $\p^{l^c}$ is in $W(k,E_{k,\mu,\tau})$, which is equal to $W(k,l^{n-1})$ by Lemma \ref{Eciclo}.
\end{proof}

\begin{lemma}\label{Wram}
Let $K/k$ be a tamely ramified abelian extension of number fields and let $\p$ be a prime ideal in $k$ whose ramification index in $K/k$ is $e$, then $\norm_{k/\Q}(\p)\in P_\Q^\mathfrak m$, where $\mathfrak m=e\cdot p_\infty$. In particular, by Lemma 1.12 of \cite{Cobbe}, $\p\in H_{k(\zeta_e)/k}^\mathfrak m$ and so its class is in $W(k,e)$.
\end{lemma}
\begin{proof}
This is Lemma I.2.1 of \cite{Endo}.
\end{proof}

\begin{lemma}\label{ln1}
Let $K/k$ be a tame $C(l^n)\rtimes_{\mu} C(l)$-extension of number fields and let $\p$ be a ramifying prime, with ramification index $e_\p$. Then the class of
\[\p^{\frac{e_\p-1}{2}\frac{l^{n+1}}{e_\p}}\]
and the class of
\[\p^{\frac{l-1}{2}\frac{l^{n+1}}{e_\p}}\]
are both in
\[W(k,l^{n-1})^{\frac{l-1}{2}l}.\]
\end{lemma}

\begin{proof}
The Galois group of $K/k$ is $C(l^n)\rtimes_{\mu} C(l)$, i.e.
\[G=\langle\sigma,\tau:\ \sigma^l=\tau^{l^n}=1,\sigma\tau\sigma^{-1}=\tau^{l^{n-1}+1}\rangle.\]

Since the ramification is tame, the inertia group at $\p$ is cyclic, generated by an element $\tau^a\sigma^b$; by induction we obtain
\[(\tau^a\sigma^b)^m=\tau^{am+abl^{n-1}(m-1)m/2}\sigma^{bm}.\]
The order $e_\p$ of $\tau^a\sigma^b$ must be a multiple of $l$, since the element $\tau^a\sigma^b$ is nontrivial and $G$ is an $l$-group. Hence, recalling that $\tau^{l^n}=1$, we obtain that $e_\p$ is the smallest positive integer such that
\[\tau^{ae_\p}\sigma^{be_\p}=1.\]
First of all we assume that $l^2$ divides $e_\p$. If $l^\beta$ is the exact power of $l$ dividing $a$, we obtain that $e_\p=l^{n-\beta}$ and in particular that $\beta\leq n-2$. So we have
\[\sigma_*(\tau^a\sigma^b)=\sigma\tau^a\sigma^b\sigma^{-1}=\tau^{a(l^{n-1}+1)}\sigma^b=(\tau^a\sigma^b)^{l^{n-1}+1}\]
and
\[\tau_*(\tau^a\sigma^b)=\tau\tau^a\sigma^b\tau^{-1}=\tau^{a-bl^{n-1}}\sigma^b=(\tau^a\sigma^b)^{-\tilde abl^{n-1-\beta}+1},\]
where $a\tilde a\equiv l^\beta\pmod{l^n}$. Hence, in particular, the inertia group is a normal subgroup of $G$. Thus we can decompose our extension in $K/k_1$ and $k_1/k$ which are both Galois and such that $\p$ is totally ramified in $K/k_1$ and unramified in $k_1/k$. By Lemma 2.14 of \cite{Cobbe} the class of $\p$ is in $W(k,E_{k,\rho,\tau^a\sigma^b})$, where the action $\rho$ is induced by the conjugation in $G$ and, in particular, it sends $\tau$ to the elevation to the $-\tilde abl^{n-1-\beta}+1$-th power, as seen above, and $\sigma$ to the elevation to the $l^{n-1}+1$-th power. The group $G_{k,\rho,\tau^a\sigma^b}$ consists of those elements $g$ of $\gal(k(\zeta_{l^{n-\beta}})/k)$ such that $\nu_{k,\tau^a\sigma^b}(g)$ is congruent to a product of powers of $l^{n-1}+1$ and $-\tilde abl^{n-1-\beta}+1$ modulo $l^{n-\beta}$. But all these are congruent to $1$ modulo $l^{n-1-\beta}$ and thus $G_{k,\rho,\tau^a\sigma^b}|_{k(\zeta_{l^{n-1-\beta}})}=\{1\}$. Hence
\[E_{k,\rho,\tau^a\sigma^b}\supseteq k(\zeta_{l^{n-1-\beta}})\supseteq k\left(\zeta_\frac{e_\p}{l}\right)\]
i.e.
\[W(k,E_{k,\rho,\tau^a\sigma^b})\subseteq W\left(k,\frac{e_\p}{l}\right).\]
Hence, by the assumption that $l^2|e_\p$, the class of 
\[\p^{\frac{l-1}{2}\frac{l^{n+1}}{e_\p}}\]
is in
\[W\left(k,\frac{e_\p}{l}\right)^{\frac{l-1}{2}\frac{l^{n+1}}{e_\p}}\subseteq W\left(k,l^{n-1}\right)^{\frac{l-1}{2}l}\]
and the same is true for
\[\p^{\frac{e_\p-1}{2}\frac{l^{n+1}}{e_\p}}.\]

It remains to consider the case $e_\p=l$. We now define $k_1$ as the fixed field of $\tau$ and we first assume that $\p$ ramifies in $K/k_1$. Then its inertia group in $\gal(K/k_1)=C(l^n)$ is of order $l$ and thus must be generated by $\tau^{l^{n-1}}$. Hence by Lemma 2.14 of \cite{Cobbe} the class of $\p$ is in $W(k,E_{k,\mu,\tau^{l^{n-1}}})$ and $\p^{(l-1)\frac{l^{n+1}}{e_\p}}$ is the square of an ideal in $W(k,E_{k,\mu,\tau^{l^{n-1}}})^{\frac{l-1}{2}l^n}$, which is contained in $W(k,E_{k,\mu,\tau})^{\frac{l-1}{2}l}$ by Lemma \ref{troppo}. Hence, by Lemma \ref{Eciclo}, the class of
\[\p^{\frac{l-1}{2}\frac{l^{n+1}}{e_\p}}=\p^{\frac{e_\p-1}{2}\frac{l^{n+1}}{e_\p}}\]
is in \[W(k,l^{n-1})^{\frac{l-1}{2}l}.\]

Finally let us consider the case of $\p$ ramified in $k_1/k$. By Lemma \ref{Wram} the class of $\p$ is in $W(k,l)$. Hence the class of
\[\p^{\frac{l-1}{2}\frac{l^{n+1}}{e_\p}}=\p^{\frac{e_\p-1}{2}\frac{l^{n+1}}{e_\p}}\]
is in
\[W(k,l)^{\frac{l-1}{2}l^{n}}.\]
By Lemma \ref{Wexp}
\[W(k,l)^{\frac{l-1}{2}l^{n}}\subseteq W(k,l^{n-1})^{\frac{l-1}{2}l^{2}}\subseteq W(k,l^{n-1})^{\frac{l-1}{2}l}.\]
\end{proof}

\begin{teo}\label{main}
We have
\[\rt(k,C(l^n)\rtimes_{\mu} C(l))= W(k,l^{n-1})^{\frac{l-1}{2}l}.\]
Further the group $C(l^n)\rtimes_{\mu} C(l)$ is good.
\end{teo}

\begin{proof}
By Theorems \ref{discriminante} and \ref{stdisc}, by Lemma \ref{lncostruttivo} and Lemma \ref{ln1} it is immediate that
\[\rt(k,C(l^n)\rtimes_{\mu} C(l))= W(k,l^{n-1})^{\frac{l-1}{2}l}.\]
The prove that $C(l^n)\rtimes_{\mu} C(l)$ is good is now straightforward using the same results.
\end{proof}
\end{section}

\begin{section}{Nonabelian extensions of order $l^3$}
As a particular case of Theorem \ref{main} we state the following proposition.
\begin{prop}
The group $C(l^2)\rtimes_{\mu} C(l)$ is good and
\[\rt(k,C(l^2)\rtimes_{\mu} C(l))= W(k,l)^{\frac{l-1}{2}l}.\]
\end{prop}
Up to isomorphism, the only other nonabelian group of order $l^3$ is
\[G=\langle x,y,\sigma:\ x^l=y^l=\sigma^l=1,\sigma x=x\sigma, \sigma y=y\sigma, yx=xy\sigma\rangle,\]
which is a semidirect product of the normal subgroup $\langle x,\sigma\rangle\cong C(l)\times C(l)$ and the cyclic subgroup $\langle y \rangle$ of order $l$, where the action $\mu_1$ is given by conjugation.
Cl\'ement Bruche proved in \cite{Bruche} that
\[\rt(k,G)=W(k,l)^{\frac{l-1}{2}l^2}.\]

We can give a different proof of Bruche's result, using class field theory. We will also prove that the nonabelian group of order $l^3$ and exponent $l$ studied by Bruche is a good group.
\begin{lemma}\label{costruttivoBruche}
Let $k$ be a number field, then
\[\rt(k,G)\supseteq W(k,l)^{\frac{l-1}{2}l^2}.\]
Further, for any $x\in W(k,l)$ and any positive integer $a$, there exists a tame $G$-extension of $k$ with Steinitz class $x^{\frac{l-1}{2}l^2}$ and such that any nontrivial subextension of $K/k$ is ramified at some primes which are unramified in $k(\zeta_a)/k$.
\end{lemma}

\begin{proof}
Let $x\in W(k,l)$ and $n\in\N\setminus\{0\}$. By Theorem 2.19 in \cite{Cobbe} there exists a $C(l)$-extension $k_1$ with Steinitz class $x^{l-1}$ and which is totally ramified at some prime ideals, which are unramified in $k(\zeta_a)/k$. Let $\p$ be one of them.

Now we would like to use Lemma 2.10 of \cite{Cobbe} to obtain a $C(l)\times C(l)$ extension of $K/k_1$ which is Galois over $k$, with $\gal(K/k)\cong G$. Unfortunately this is not possible since the exact sequence
\[1\to C(l)\times C(l)\to \mathcal H\to C(l)\to 1\]
does not identify the group $\mathcal H$ uniquely as the group $G$. Nevertheless, the argument of that lemma
at least produces a $C(l)\times C(l)$-extension of $k_1$, which is Galois over $k$ and with $\st(K/k_1)=1$. Further we get that $\gal(K/k)$ is nonabelian of oder $l^3$ (since
the action of $C(l)$ on $C(l)\times C(l)$ is the given one and in particular not trivial), that $K/k_1$ is unramified at $\p$ and that any nontrivial subextension of $K/k$ is ramified at some primes which are unramified in $k(\zeta_a)/k$.

We want to prove that $\gal(K/k)\cong G$. To this aim, we assume that this is not the case, i.e. that $\gal(K/k)\cong C(l^2)\rtimes_\mu C(l)$, and we derive a contradiction. First of all, by construction, $\gal(K/k_1)\cong C(l)\times C(l)$ and this must be a subgroup of $\gal(K/k)\cong C(l^2)\rtimes_\mu C(l)$: the only possibility is that it is the
subgroup $H$ which arises by replacing $C(l^2)$ (the left hand factor in the semidirect product) by its subgroup of order $l$; H happens to consist of all elements of $C(l^2)\rtimes_\mu C(l)$ having order $1$ or $l$. Since the prime ideal $\p$ ramifies in $k_1/k$ and not in $K/k_1$, its ramification index is $l$ and, therefore, its inertia group is contained in $H$. Hence by Galois theory we conclude that the inertia field of $\p$ in $K/k$ contains $k_1$, i.e. that $\p$ ramifies in $K/k_1$ and not in $k_1/k$. This is a contradiction, since $\p$ is ramified in $k_1/k$.

Hence we have proved that in the above construction the extension $K/k$ has Galois group $G$. By Proposition \ref{stintermediateextension},
\[\st(K/k)=\st(k_1/k)^{[K:k_1]}=x^{\frac{l-1}{2}l^2}\norm_{k_1/k}(\st(K/k_1)).\]
\end{proof}

To prove the opposite inclusion we need the following lemma.
\begin{lemma}\label{ramifl}
Let $K/k$ be a tame $G$-extension of number fields. The ramification index of a prime ramifying in $K/k$ is $l$ and its class is contained in $W(k,l)$.
\end{lemma}
\begin{proof}
The ramification index of a ramifying prime is equal to $l$, since the corresponding inertia group must be cyclic and any nontrivial element in $G$ is of order $l$.

Let $k_1$ be the subfield of $K$ fixed by the normal abelian subgroup $\langle x,\sigma\rangle$ of the Galois group $G$ of $K/k$.

If a prime $\p$ ramifies in $k_1/k$, then its class is in $W(k,l)$ by Lemma \ref{Wram}.

If a prime $\p$ ramifies in $K/k_1$, then it is unramified in $k_1/k$ (the ramification index is prime) and so its inertia group is generated by an element of the form $x^a\sigma^c$, where $a,c\in\{0,1,\dots,l-1\}$ are not both $0$. By Lemma 2.14 of \cite{Cobbe} the class of $\p$ is in $W(k,E_{k,\mu_1,x^a\sigma^c})$. For any $b\in\{0,1,\dots,l-1\}$ we have
\[\mu_1(y^b)(x^a\sigma^c)=y^bx^a\sigma^c y^{-b}=x^a\sigma^{c+ab}\]
and this expression cannot be a nontrivial power of $a^a\sigma^c$. Hence, by definition, the group $G_{k,\mu_1,x^a\sigma^c}$ must be trivial and we conclude that $E_{k,\mu_1,x^a\sigma^c}=k(\zeta_l)$. Therefore, in particular, the class of the prime ideal $\p$ is contained in $W(k,l)$.
\end{proof}

\begin{prop}
The group $G$ is good and
\[\rt(k,G)= W(k,l)^{\frac{l-1}{2}l^2}.\]
\end{prop}
\begin{proof}
One inclusion is given by Lemma \ref{costruttivoBruche}. The proof that $\rt(k,G)\subseteq W(k,l)^{\frac{l-1}{2}l^2}$ follows by Lemma \ref{ramifl}, since for any tame $G$-extension $K/k$ of number fields the Steinitz class is the class of the ideal
\[\prod_{\p:\ e_\p\neq 1} \p^{\frac{e_\p-1}{2}\frac{l^3}{e_\p}}=\prod_{\p:\ e_\p\neq 1} \p^{\frac{l-1}{2}l^2},\]
which is contained in $W(k,l)^{\frac{l-1}{2}l^2}$. Now we prove that all the properties of good groups are verified.
\begin{enumerate}
\item This is clear, since $W(k,l)$ is a group.
\item For any prime $\p$, ramifying in a tame $G$-extension $K/k$ of number fields, by Lemma \ref{ramifl}, the class of
\[\p^{\frac{e_\p-1}{2}\frac{l^3}{e_\p}}=\p^{\frac{l-1}{2}l^2}\]
is contained in $W(k,l)^{\frac{l-1}{2}l^2}$, which is equal to $\rt(k,G)$.
\item Immediate by Lemma \ref{ramifl} and the explicit formula for $\rt(k,G)$.
\item This follows by Lemma \ref{costruttivoBruche}.
\end{enumerate}
\end{proof}
\end{section}

\nocite{McCulloh}
\nocite{Long0}
\nocite{Long2}
\nocite{Endo}
\nocite{Carter}
\nocite{Massy}
\nocite{Sodaigui1}
\nocite{Sodaigui2}
\nocite{Soverchia}
\nocite{GodinSodaigui_A4}
\nocite{GodinSodaigui_ottaedri}
\nocite{ByottGreitherSodaigui} 
\nocite{CarterSodaigui_quaternionigeneralizzati}
\nocite{BrucheSodaigui}

\bibliography{bibsteinitz}

\begin{thebibliography}{10}

\bibitem{Bruche}
C.~Bruche.
\newblock Classes de {S}teinitz d'extensions non ab\'eliennes de degr\'e
  {$p^3$}.
\newblock {\em Acta Arith.}, 137(2):177--191, 2009.

\bibitem{BrucheSodaigui}
C.~Bruche and B.~Soda{\"{\i}}gui.
\newblock On realizable {G}alois module classes and {S}teinitz classes of
  nonabelian extensions.
\newblock {\em J. Number Theory}, 128(4):954--978, 2008.

\bibitem{ByottGreitherSodaigui}
N.~P. Byott, C.~Greither, and B.~Soda{\"{\i}}gui.
\newblock Classes r\'ealisables d'extensions non ab\'eliennes.
\newblock {\em J. Reine Angew. Math.}, 601:1--27, 2006.

\bibitem{Carter}
J.~E. Carter.
\newblock Steinitz classes of a nonabelian extension of degree {$p\sp 3$}.
\newblock {\em Colloq. Math.}, 71(2):297--303, 1996.

\bibitem{Carter2}
J.~E. Carter.
\newblock Steinitz classes of nonabelian extensions of degree {$p^3$}.
\newblock {\em Acta Arith.}, 78(3):297--303, 1997.

\bibitem{CarterSodaigui_quaternionigeneralizzati}
J.~E. Carter and B.~Soda{\"{\i}}gui.
\newblock Classes de {S}teinitz d'extensions quaternioniennes
  g\'en\'eralis\'ees de degr\'e {$4p\sp r$}.
\newblock {\em J. Lond. Math. Soc. (2)}, 76(2):331--344, 2007.

\bibitem{tesi}
A.~Cobbe.
\newblock {\em Steinitz classes of tamely ramified Galois extensions of
  algebraic number fields}.
\newblock PhD thesis, Scuola Normale Superiore, Pisa, 2010.

\bibitem{Cobbe}
A.~Cobbe.
\newblock {\em Steinitz classes of tamely ramified Galois extensions of
  algebraic number fields}.
\newblock arXiv:0910.5080v1, to appear in Journal of Number Theory.

\bibitem{Endo}
L.~P. Endo.
\newblock {\em Steinitz classes of tamely ramified Galois extensions of
  algebraic number fields}.
\newblock PhD thesis, University of Illinois at Urbana-Champaign, 1975.

\bibitem{GodinSodaigui_A4}
M.~Godin and B.~Soda{\"{\i}}gui.
\newblock Classes de {S}teinitz d'extensions \`a groupe de {G}alois {$A\sb 4$}.
\newblock {\em J. Th\'eor. Nombres Bordeaux}, 14(1):241--248, 2002.

\bibitem{GodinSodaigui_ottaedri}
M.~Godin and B.~Soda{\"{\i}}gui.
\newblock Module structure of rings of integers in octahedral extensions.
\newblock {\em Acta Arith.}, 109(4):321--327, 2003.

\bibitem{Lang}
S.~Lang.
\newblock {\em Algebraic number theory}.
\newblock GTM 110. Springer-Verlag, New York, second edition, 1994.

\bibitem{Long2}
R.~Long.
\newblock Steinitz classes of cyclic extensions of degree {$l\sp{r}$}.
\newblock {\em Proc. Amer. Math. Soc.}, 49:297--304, 1975.

\bibitem{Long0}
R.~L. Long.
\newblock Steinitz classes of cyclic extensions of prime degree.
\newblock {\em J. Reine Angew. Math.}, 250:87--98, 1971.

\bibitem{Massy}
R.~Massy and B.~Soda{\"{\i}}gui.
\newblock Classes de {S}teinitz et extensions quaternioniennes.
\newblock {\em Proyecciones}, 16(1):1--13, 1997.

\bibitem{McCulloh}
L.~R. McCulloh.
\newblock Cyclic extensions without relative integral bases.
\newblock {\em Proc. Amer. Math. Soc.}, 17:1191--1194, 1966.

\bibitem{McCulloh_Crelle}
L.~R. McCulloh.
\newblock Galois module structure of abelian extensions.
\newblock {\em J. Reine Angew. Math.}, 375/376:259--306, 1987.

\bibitem{Narkiewicz}
W.~Narkiewicz.
\newblock {\em Elementary and analytic theory of algebraic numbers}.
\newblock Springer Monographs in Mathematics. Springer-Verlag, Berlin, third
  edition, 2004.

\bibitem{Neukirch}
J.~Neukirch.
\newblock {\em Class field theory}, volume 280 of {\em Grundlehren der
  Mathematischen Wissenschaften [Fundamental Principles of Mathematical
  Sciences]}.
\newblock Springer-Verlag, Berlin, 1986.

\bibitem{Sodaigui1}
B.~Soda{\"{\i}}gui.
\newblock Classes de {S}teinitz d'extensions galoisiennes relatives de degr\'e
  une puissance de 2 et probl\`eme de plongement.
\newblock {\em Illinois J. Math.}, 43(1):47--60, 1999.

\bibitem{Sodaigui2}
B.~Soda{\"{\i}}gui.
\newblock Relative {G}alois module structure and {S}teinitz classes of dihedral
  extensions of degree {$8$}.
\newblock {\em J. Algebra}, 223(1):367--378, 2000.

\bibitem{Soverchia}
E.~Soverchia.
\newblock Steinitz classes of metacyclic extensions.
\newblock {\em J. London Math. Soc. (2)}, 66(1):61--72, 2002.

\end{thebibliography}
\addcontentsline{toc}{section}{Bibliography}
\bibliographystyle{abbrv}

\end{document}